\documentclass[brochure,english,12pt]{bourbaki}
\usepackage[matrix,arrow]{xy}
\usepackage{amssymb,amsfonts,amsmath}
\usepackage{euscript}   
\usepackage{eepic}
\usepackage[francais]{babel}
\DeclareMathAlphabet{\eurm}{U}{eur}{m}{n}

\numberwithin{equation}{section}

\newcommand{\bean}{\begin{eqnarray}}
\newcommand{\eean}{\end{eqnarray}}
\newcommand{\be}{\begin{displaymath}}
\newcommand{\ee}{\end{displaymath}}
\newcommand{\bea}{\begin{eqnarray*}}
\newcommand{\eea}{\end{eqnarray*}}

\newcommand{\secref}[1]{Section~\ref{#1}}

\newcommand{\nc}{\newcommand}
\nc{\on}{\operatorname}
\nc{\ch}{\mbox{ch}}
\nc{\Z}{{\mathbb Z}}
\nc{\C}{{\mathbb C}}
\nc{\pone}{{\mathbb P}^1}
\nc{\pa}{\partial}
\nc{\F}{{\mathcal F}}
\nc{\arr}{\rightarrow}
\nc{\larr}{\longrightarrow}
\nc{\al}{\alpha}
\nc{\ri}{\rangle}
\nc{\lef}{\langle}
\nc{\W}{{\mathcal W}}
\nc{\la}{\lambda}
\nc{\ep}{\epsilon}
\nc{\su}{\widehat{{\mathfrak s}{\mathfrak l}}_2}
\nc{\sw}{{\mathfrak s}{\mathfrak l}}
\nc{\g}{{\mathfrak g}}
\nc{\h}{{\mathfrak h}}
\nc{\n}{{\mathfrak n}}
\nc{\N}{\widehat{\n}}
\nc{\G}{\widehat{\g}}
\nc{\De}{\Delta}
\nc{\gt}{\widetilde{\g}}
\nc{\Ga}{\Gamma}
\nc{\one}{{\mathbf 1}}
\nc{\z}{{\mathfrak Z}}
\nc{\La}{\Lambda}
\nc{\wt}{\widetilde}
\nc{\wh}{\widehat}
\nc{\cri}{_{\kappa_c}}
\nc{\kk}{h^\vee}
\nc{\sun}{\widehat{\sw}_N}
\nc{\si}{\sigma}
\nc{\el}{\ell}
\nc{\bi}{\bibitem}
\nc{\om}{\omega}
\nc{\ol}{\overline}
\nc{\ds}{\displaystyle}
\nc{\dzz}{\frac{dz}{z}}
\nc{\Res}{\on{Res}}
\nc{\mc}{\mathcal}
\nc{\Cal}{\mathcal}
\nc{\bb}{{\mathfrak b}}
\nc{\ot}{\otimes}
\nc{\R}{{\mc R}}
\nc{\yy}{{\mc Y}}
\nc{\ga}{\gamma}

\nc{\us}{\underset}
\nc{\opl}{\oplus}
\nc{\beq}{\begin{equation}}
\nc{\Fq}{{\mathbb F}_q}
\nc{\Mq}{{\mathcal M}}
\nc{\Rep}{\on{Rep}}
\nc{\sssec}{\subsubsection}
\nc{\ssec}{\subsection}
\nc{\lan}{\langle}
\nc{\ran}{\rangle}

\nc{\D}{\mathcal D} \nc{\Vect}{\on{Vect}} \nc{\ghat}{\G}
\nc{\T}{\mc T} \nc{\Tloc}{\T^\g_{\on{loc}}} \nc{\vac}{|0\ran}
\nc{\Wick}{{\mb :}} \nc{\mb}{\mathbf} \nc{\delz}{\partial_z}
\nc{\K}{{\cali K}} \nc{\cali}{\mathcal} \nc{\li}{\mathfrak l}
\nc{\lt}{\widetilde{\li}} \nc{\astar}{a^*} \nc{\cA}{{\mc A}}
\nc{\ka}{\kappa}

\def\FF{\text{\bf \em F}}  
\def\PP{\text{\bf \em P}}  

\nc{\OO}{{\mc O}}
\nc{\AutO}{\on{Aut}\OO}
\nc{\DerO}{\on{Der}\OO}
\nc{\DerpO}{\on{Der}_+\OO}
\nc{\Au}{{\mc A}ut}
\nc{\mf}{\mathfrak}
\nc{\V}{{\mathbb V}}
\nc{\hh}{\wh{\h}}

\nc{\pp}{{\mathfrak p}}
\nc{\mm}{{\mathfrak m}}
\nc{\rr}{{\mathfrak r}}
\nc{\ket}{\rangle}
\nc{\zz}{{\mathfrak z}}
\nc{\gr}{\on{gr}}
\nc{\Spe}{\on{Spec}}
\nc{\rv}{\crho}
\nc{\can}{\on{can}}
\nc{\Op}{\on{Op}_G(D)}
\nc{\MOp}{\on{MOp}_G(D)}
\nc{\Db}{{\mathbb D}}
\nc{\ww}{w}

\nc{\af}{{\mathbb A}^1}
\nc{\bs}{\backslash}
\nc{\laa}{(\la_i)}
\nc{\zn}{(z_i)}

\nc{\cla}{\check{\la}}
\nc{\cmu}{\check{\mu}}
\nc{\crho}{\check{\rho}}
\nc{\chal}{\check{\al}}
\nc{\cc}{{\mathfrak c}}

\nc{\M}{{\mathbb M}}

\nc{\ZZ}{{\mc Z}}

\nc{\UU}{{\mathbb U}}

\nc{\Conn}{\on{Conn}(\Omega^{\crho})}
\nc{\Con}{\on{Conn}(\Omega^{-\rho})}
\nc{\Co}{\on{Conn}(\Omega^{\rho})}

\nc{\ppart}{(\!(t)\!)}
\nc{\pparl}{(\!(\la)\!)}
\nc{\zpart}{(\!(z)\!)}
\nc{\ppzi}{(\!(t-z_i)\!)}
\nc{\ppinf}{(\!(t^{-1})\!)}
\nc{\Ind}{\on{Ind}}
\nc{\I}{{\mathbb I}}

\nc{\Bun}{\on{Bun}}

\nc{\CC}{C}

\nc{\gtil}{\wt{\g}}
\nc{\ntil}{\wt{\n}}
\nc{\htil}{\wt{\h}}
\nc{\gbar}{\ol{\g}}
\nc{\nbar}{\ol{\n}}
\nc{\bbar}{\ol{\bb}}

\def\M{{\mathcal M}}
\def\MH{{\mathcal M}_H}

\def\L{{\mathcal L}}

\def\R{{\mb R}}
\def\cal{\mathcal}

\def\CO{{\cal O}}
\def\O{\CO}

\nc{\mbb}{\mathbb}
\def\neg{\negthinspace}
\nc{\LG}{{}^L\neg G}
\nc{\LH}{{}^L\neg H}
\nc{\LZ}{{}^L\neg Z}

\def\B{{\cal B}}
\def\D{{\cal D}}

\def\A{{\cal A}}
\def\N{{\cal N}}

\def\M{{\cal M}}
\def\L{{\cal L}}
\nc{\Pic}{\on{Pic}}

\nc{\Irrep}{\on{Irrep}}

\def\E{{\mathcal E}}

\def\R{{\cal R}}
\def\T{{\cal T}}

\def\Pic{{\rm Pic}}

\nc{\Q}{{\mathbb Q}}
\nc{\Qbar}{\ol\Q}
\nc{\Ql}{{\mathbb Q}_\ell}
\nc{\Gal}{\on{Gal}}
\nc{\AD}{{\mathbb A}}
\nc{\hl}{h^{\leftarrow}}
\nc{\hr}{h^{\rightarrow}}
\nc{\supp}{\on{supp}}
\nc{\Loc}{\on{Loc}}

\date{Juin 2009}
\bbkannee{61\`eme ann\'ee, 2008-2009}
\bbknumero{1010}
\title{GAUGE THEORY AND LANGLANDS DUALITY}
\author{Edward FRENKEL}\thanks{Supported by DARPA through the
  grant HR0011-09-1-0015 and by Fondation Sciences Math\'ematiques de
  Paris}
\address{Department of Mathematics\\
University of California\\
Berkeley, CA 94720-3840 -- U.S.A.}
\email{frenkel@math.berkeley.edu}

\begin{document}
\maketitle

\noindent{\bf INTRODUCTION}

\bigskip

In the late 1960s Robert Langlands launched what has become known as
the Langlands Program with the ambitious goal of relating deep
questions in Number Theory to Harmonic Analysis \cite{L}. In
particular, Langlands conjectured that Galois representations and
motives can be described in terms of the more tangible data of
automorphic representations. A striking application of this general
principle is the celebrated Shimura--Taniyama--Weil conjecture (which
implies Fermat's Last Theorem), proved by A. Wiles and others, which
says that information about Galois representations associated to
elliptic curves over ${\mathbb Q}$ is encoded in the Fourier expansion
of certain modular forms on the upper-half plane.

One of the most fascinating and mysterious aspects of the Langlands
Program is the appearance of the {\em Langlands dual group}. Given a
reductive algebraic group $G$, one constructs its Langlands dual $\LG$
by applying an involution to its root data. Under the Langlands
correspondence, automorphic representations of the group $G$
correspond to Galois representations with values in $\LG$.

Surprisingly, the Langlands dual group also appears in Quantum Physics
in what looks like an entirely different context; namely, the {\em
electro-magnetic duality}. Looking at the Maxwell equations describing
the classical electromagnetism, one quickly notices that they are
invariant under the exchange of the electric and magnetic fields. It
is natural to ask whether this duality exists at the quantum level. In
quantum theory there is an important parameter, the electric charge
$e$. Physicists have speculated that there is an electro-magnetic
duality in the quantum theory under which $e \longleftrightarrow
1/e$. Under this duality the electrically charged particle should be
exchanged with a magnetically charged particle, called magnetic
monopole, first theorized by P. Dirac (so far, it has not been
discovered experimentally).

In modern terms, Maxwell theory is an example of 4D {\em gauge theory}
(or Yang--Mills theory) which is defined, classically, on the space of
connections on various $G_c$-bundles on a four-manifold $M$, where
$G_c$ is a compact Lie group.\footnote{We will use the notation $G$
for a complex Lie group and $G_c$ for its compact form. Note that
physicists usually denote by $G$ a compact Lie group and by $G_{\C}$
its complexification.}  Electromagnetism corresponds to the simplest,
abelian, compact Lie group $U(1)$. It is natural to ask whether there
is a non-abelian analogue of the electro-magnetic duality for gauge
theories with non-abelian gauge groups.

The answer was proposed in the late 1970s, by Montonen and Olive
\cite{MO}, following Goddard, Nuyts and Olive \cite{GNO} (see also
\cite{EW,O}). A gauge theory has a coupling constant $g$, which plays
the role of the electric charge $e$. The conjectural non-abelian
electro-magnetic duality, which has later become known as $S$-{\em
duality}, has the form
\begin{equation}    \label{S-duality}
(G_c,g) \longleftrightarrow ({}\LG_c,1/g).
\end{equation}
In other words, the duality states that the gauge theory with gauge
group $G_c$ (more precisely, its ``$N=4$ supersymmetric'' version) and
coupling constant $g$ should be equivalent to the gauge theory with
the Langlands dual gauge group $\LG_c$ and coupling constant $1/g$
(note that if $G_c=U(1)$, then $\LG_c$ is also $U(1)$). If true, this
duality would have tremendous consequences for quantum gauge theory,
because it would relate a theory at small values of the coupling
constant (weak coupling) to a theory with large values of the coupling
constant (strong coupling). Quantum gauge theory is usually defined as
a power series expansion in $g$, which can only converge for small
values of $g$. It is a very hard problem to show that these series
make sense beyond perturbation theory. $S$-duality indicates that the
theory does exist non-perturbatively and gives us a tool for
understanding it at strong coupling. That is why it has become a holy
grail of modern Quantum Field Theory.

Looking at \eqref{S-duality}, we see that the Langlands dual group
shows up again. Could it be that the Langlands duality in Mathematics
is somehow related to $S$-duality in Physics?

This question has remained a mystery until about five years ago. In
March of 2004, DARPA sponsored a meeting of a small group of
physicists and mathematicians at the Institute for Advanced Study in
Princeton (which I co-organized) to tackle this question. At the end
of this meeting Edward Witten gave a broad outline of a relation
between the two topics. This was explained in more detail in his
subsequent joint work \cite{KW} with Anton Kapustin. This paper, and
the work that followed it, opened new bridges between areas of great
interest for both physicists and mathematicians, leading to new ideas,
insights and directions of research.

The goal of these notes is to describe briefly some elements of the
emerging picture. In Sections \ref{LP} and \ref{GLP}, we will discuss
the Langlands Program and its three flavors, putting it in the context
of Andr\'e Weil's ``big picture''. This will eventually lead us to a
formulation of the geometric Langlands correspondence as an
equivalence of certain categories of sheaves in \secref{cat ver}. In
\secref{PHY} we will turn to the $S$-duality in topological twisted
$N=4$ super--Yang--Mills theory. Its dimensional reduction gives rise
to the Mirror Symmetry of two-dimensional sigma models associated to
the Hitchin moduli spaces of Higgs bundles. In \secref{MIRROR} we will
describe a connection between the geometric Langlands correspondence
and this Mirror Symmetry, following \cite{KW}, as well as its ramified
analogue \cite{GW}. In \secref{BRANES} we will discuss subsequent work
and open questions.

\medskip

\noindent {\bf Acknowledgments.} I thank Sergei Gukov, Vincent
Lafforgue, Robert Langlands, and Edward Witten for inspiring
discussions. I also thank S. Gukov and V. Lafforgue for their comments
on a draft of this paper.

I am grateful to DARPA (and especially Benjamin Mann) for generous
support which has been instrumental not only for my research, but
for the development of this whole area. I also thank Fondation
Sciences Math\'ematiques de Paris for its support during my stay in
Paris.

\section{LANGLANDS PROGRAM}    \label{LP}

In 1940 Andr\'e Weil was put in jail for his refusal to serve in the
army. There, he wrote a letter to his sister Simone Weil (a noted
philosopher) in response to her question as to what really interested
him in his work \cite{Weil}. This is a remarkable document, in which
Weil tries to explain, in fairly elementary terms (presumably,
accessible even to a philosopher), the ``big picture'' of mathematics,
the way he saw it. I think this sets a great example to follow for all
of us.

Weil writes about the role of {\em analogy} in mathematics, and he
illustrates it by the analogy that interested him the most: between
Number Theory and Geometry.

On one side we look at the field $\Q$ of rational numbers and its
algebraic closure $\Qbar$, obtained by adjoining all roots of all
polynomial equations in one variables with rational coefficients (like
$x^2+1=0$). The group of field automorphisms of $\Qbar$ is the {\em
Galois group} $\Gal(\Qbar/\Q)$. We are interested in the structure of
this group and its finite-dimensional representations. We may also
take a more general {\em number field} -- that is, a finite extension
$F$ of $\Q$ (such as $\Q(i)$) -- and study its Galois group and its
representations.

On the other side we have Riemann surfaces: smooth compact orientable
surfaces equipped with a complex structure, and various geometric
objects associated to them: vector bundles, their endomorphisms,
connections, etc.

At first glance, the two subjects are far apart. However, it turns out
that there are many analogies between them. The key point is that
there is another class of objects which are in-between the two. A
Riemann surface may be viewed as the set of points of a projective
algebraic curve over $\C$. In other words, Riemann surfaces may be
described by algebraic equations, such as the equation
\begin{equation}    \label{ell curve}
y^2 = x^3 + ax + b,
\end{equation}
where $a,b \in \C$. The set of complex solutions of this equation (for
generic $a,b$ for which the polynomial on the right hand side has no
multiple roots), compactified by a point at infinity, is a Riemann
surface of genus $1$. However, we may look at the equation \eqref{ell
curve} not only over $\C$, but also over other fields -- for
instance, over finite fields.

Recall that there is a unique, up to an isomorphism, finite field
$\Fq$ of $q$ elements for all $q$ of the form $p^n$, where $p$ is a
prime. In particular, ${\mathbb F}_p = \Z/p\Z \simeq \{ 0,1,\ldots,p-1
\}$, with the usual arithmetic modulo $p$. Let $a,b$ be elements of
$\Fq$. Then the equation \eqref{ell curve} defines a curve over
$\Fq$. These objects are clearly analogous to algebraic curves over
$\C$ (that is, Riemann surfaces). But there is also a deep analogy
with number fields!

Indeed, let $X$ be a curve over $\Fq$ (such as an elliptic curve
defined by \eqref{ell curve}) and $F$ the field of rational functions
on $X$. This {\em function field} is very similar to a number
field. For instance, if $X$ is the projective line over $\Fq$, then
$F$ consists of all fractions $P(t)/Q(t)$, where $P$ and $Q$ are two
relatively prime polynomials in one variable with coefficients in
$\Fq$. The ring $\Fq[t]$ of polynomials in one variable over $\Fq$ is
similar to the ring of integers and so the fractions $P(t)/Q(t)$ are
similar to the fractions $p/q$, where $p,q \in \Z$.

Thus, we find a {\em bridge}, or a ``turntable'' -- as Weil calls it --
between Number Theory and Geometry, and that is the theory of
algebraic curves over finite fields.

In other words, we can talk about three parallel tracks

$$
\begin{matrix}
\text{{\em Number Theory}} \qquad & \qquad \text{{\em Curves over}
  $\Fq$} \qquad & \qquad \text{{\em Riemann
  Surfaces}}
\end{matrix}
$$

\bigskip

Weil's idea is to exploit it in the following way: take a statement in
one of the three columns and translate it into statements in the other
columns \cite{Weil}: ``my work consists in deciphering a trilingual
text; of each of the three columns I have only disparate fragments; I
have some ideas about each of the three languages: but I know as well
there are great differences in meaning from one column to another, for
which nothing has prepared me in advance. In the several years I have
worked at it, I have found little pieces of the dictionary.'' Weil
went on to find one of the most spectacular applications of this
``Rosetta stone'': what we now call the {\em Weil conjectures}
describing analogues of the Riemann Hypothesis in Number Theory in the
context of algebraic curves over finite fields.

\medskip

It is instructive to look at the Langlands Program through the prism
of Weil's big picture. Langlands' original formulation \cite{L}
concerned the two columns on the left. Part of the Langlands Program
may be framed as the question of describing $n$-dimensional
representations of the Galois group $\on{Gal}(\ol{F}/F)$, where $F$ is
either a number field ($\Q$ or its finite extension) or the function
field of a curve over $\Fq$.\footnote{Langlands' more general
``functoriality principle'' is beyond the scope of the present
article.} Langlands proposed that such representations may be
described in terms of {\em automorphic representations} of the group
$GL_n({\mathbb A}_F)$, where $\AD_F$ is the ring of ad\`eles of $F$. I
will not attempt to explain this here referring the reader to the
surveys \cite{Gelbart,F:bull,F:houches}.

However, it is important for us to emphasize how the Langlands dual
group appears in this story. Let us replace $GL_n(\AD_F)$ by
$G(\AD_F)$, where $G$ is a general reductive algebraic group (such as
orthogonal or symplectic, or $E_8$). In the case when $G=GL_n$ its
automorphic representations are related to the $n$-dimensional
representations of $\on{Gal}(\ol{F}/F)$, that is, homomorphisms
$\on{Gal}(\ol{F}/F) \to GL_n$. The general Langlands conjectures
predict that automorphic representations of $G(\AD_F)$ are related, in
a similar way, to homomorphisms $\on{Gal}(\ol{F}/F) \to {}\LG$, where
$\LG$ is the Langlands dual group to $G$.\footnote{More precisely,
$\LG$ should be defined over $\ol{\mathbb Q}_\ell$, where $\ell$ is
relatively prime to $q$, and we should consider homomorphisms
$\on{Gal}(\ol{F}/F) \to {}\LG(\ol{\mathbb Q}_\ell)$ which are
continuous with respect to natural topology (see, e.g., Section 2.2 of
\cite{F:houches}).}

It is easiest to define $\LG$ in the case when $G$, defined over a
field $k$, is split over $k$, that is, contains a maximal split torus
$T$ (which is the product of copies of the multiplicative group $GL_1$
over $k$). We associate to $T$ two lattices: the weight lattice
$X^*(T)$ of homomorphisms $T \to GL_1$ and the coweight lattice
$X_*(T)$ of homomorphisms $GL_1 \to T$. They contain the sets of roots
$\Delta \subset X^*(T)$ and coroots $\Delta^\vee \subset X_*(T)$ of
$G$, respectively. The quadruple $(X^*(T),X_*(T),\Delta,\Delta^\vee)$
is called the root data for $G$ over $k$. The root data determines the
split group $G$ up to an isomorphism.

Let us now exchange the lattices of weights and coweights and the sets
of simple roots and coroots. Then we obtain the root data
$$
(X_*(T),X^*(T),\Delta^\vee,\Delta)
$$
of another reductive algebraic group over $\C$ (or $\ol{\mathbb
Q}_\ell$), which is denoted by $\LG$.\footnote{In Langlands'
definition \cite{L}, $\LG$ also includes the Galois group of a finite
extension of $F$. This is needed for non-split groups, but since we
focus here on the split case, this is not necessary.} Here are some
examples:

\bigskip

\begin{center}
\begin{tabular}{l|l}
$G$ & $\LG$ \\
\hline
$GL_n$ & $GL_n$ \\
\hline
$SL_n$ & $PGL_n$ \\
\hline
$Sp_{2n}$ & $SO_{2n+1}$ \\
\hline
$Spin_{2n}$ & $SO_{2n}/\Z_2$ \\
\hline
$E_8$ & $E_8$ \\
\end{tabular}
\end{center}

\bigskip


In the function field case we expect to have a correspondence between
homomorphisms\footnote{More precisely, the Galois group should be
replaced by its subgroup called the Weil group.} $\on{Gal}(\ol{F}/F)
\to {}\LG$ and automorphic representations of $G(\AD_F)$, where
$\AD_F$ is the ring of ad\`eles of $F$,
$$
\AD_F = \prod_{x \in X}{}' F_x,
$$
$F_x \simeq {\mathbb F}_{q_x}(\!(t_x)\!)$ being the completion of
the field of functions at a closed point $x$ of $X$, and the prime
means that we take the restricted product, in the sense that for all
but finitely many $x$ the element of $F_x$ belongs to its ring of
integers $\OO_x \simeq {\mathbb F}_x[[t_x]]$. We have a natural
diagonal inclusion $F \subset \AD_F$ an hence $G(F) \subset
G(\AD_F)$. Roughly speaking, an irreducible representation of
$G(\AD_F)$ is called automorphic if it occurs in the decomposition of
$L^2(G(F)\bs G(\AD_F))$ (with respect to the right action of
$G(\AD_F)$).

For $G=GL_n$, in the function field case, the Langlands correspondence
is a bijection between equivalence classes of irreducible
$n$-dimensional ($\ell$-adic) representations of $\on{Gal}(\ol{F}/F)$
(more precisely, the Weil group) and cuspidal automorphic
representations of $GL_n(\AD_F)$. It has been proved by V. Drinfeld
\cite{Dr1,Dr2} for $n=2$ and by L. Lafforgue \cite{Laf} for $n>2$. A
lot of progress has also been made recently in proving the Langlands
correspondence for $GL_n$ in the number field case.

For other groups the correspondence is expected to be much more
subtle; for instance, it is not one-to-one. Homomorphisms from the
Weil group of $F$ to $\LG$ (and more general parameters introduced by
J. Arthur, see \secref{more general}) should parametrize certain
collections of automorphic representations called ``$L$-packets.''
This has only been proved in a few cases so far.

\section{GEOMETRIC LANGLANDS CORRESPONDENCE}    \label{GLP}

The above discussion corresponds to the middle column in the Weil big
picture. What should be its analogue in the right column -- that is,
for complex curves?

In order to explain this, we need a geometric reformulation of the
Langlands correspondence which would make sense for curves defined
both over a finite field and over $\C$. Thus, we need to find
geometric analogues of the notions of Galois representations and
automorphic representations.

The former is fairly easy. Let $X$ be a curve over a field $k$ and $F
= k(X)$ the field of rational functions on $X$. If $Y \to X$ is a
covering of $X$, then the field $k(Y)$ of rational functions on $Y$ is
an extension of the field $F = k(X)$ of rational functions on $X$, and
the Galois group $\on{Gal}(k(Y)/k(X))$ may be viewed as the group of
``deck transformations'' of the cover. If our cover is {\em
unramified}, then this group is a quotient of the (arithmetic)
fundamental group of $X$. For a cover ramified at points
$x_1,\ldots,x_n$, it is a quotient of the (arithmetic) fundamental
group of $X \bs \{ x_1,\ldots,x_n \}$. From now on (with the exception
of \secref{ram}) we will focus on the unramified case. This means that
we replace $\on{Gal}(\ol{F}/F)$ by its maximal unramified quotient,
which is nothing but the (arithmetic) fundamental group of $X$. Its
geometric analogue, when $X$ is defined over $\C$, is $\pi_1(X)$.

Thus, the geometric counterpart of a (unramified) homomorphism
$\on{Gal}(\ol{F}/F) \to {}\LG$ is a homomorphism $\pi_1(X) \to {}\LG$.

{}From now on, let $X$ be a smooth projective connected algebraic
curve defined over $\C$. Let $G$ be a complex reductive algebraic
group and $\LG$ its Langlands dual group. Then homomorphisms $\pi_1(X)
\to {}\LG$ may be described in differential geometric terms as bundles
with a flat connection (the monodromy of the flat connection gives
rise to a homomorphism $\pi_1(X) \to {}\LG$). Let $E$ be a smooth
principal $\LG$-bundle on $X$. A flat connection on $E$ has two
components. The $(0,1)$ component, with respect to the complex
structure on $X$, defines holomorphic structure on $E$, and the
$(1,0)$ component defines a holomorphic connection $\nabla$. Thus, an
$\LG$-bundle with a flat connection on $X$ is the same as a pair
$(E,\nabla)$, where $E$ is a holomorphic (equivalently, algebraic)
$\LG$-bundle on $X$ and $\nabla$ is a holomorphic (equivalently,
algebraic) connection on $E$.

Thus, for complex curves the objects on one side of the Langlands
correspondence are equivalence classes of flat (holomorphic or
algebraic) $\LG$-bundles $(E,\nabla)$.

What about the other side? Here the answer is not quite as obvious. I
will sketch it briefly referring the reader to Section 3 of
\cite{F:houches} for more details.

Recall that automorphic representations of $G(\AD_F)$ (where $F$ is a
function field of a curve $X$ defined over $\Fq$) are realized in
functions on the quotient $G(F) \bs G(\AD_F)$. An unramified
automorphic representation (which corresponds to an unramified
homomorphism $\on{Gal}(\ol{F}/F) \to {}\LG$) gives rise to a function
on the double quotient $G(F) \bs G(\AD_F)/G(\OO_F)$, where $\OO_F =
\prod_{x \in X} \OO_x$. A key observation (which is due to Weil) is
that this double quotient is precisely the set of isomorphism classes
of principal $G$-bundles on our curve $X$.\footnote{From now on we
will only consider algebraic bundles.} This statement is also true if
the curve $X$ is defined over $\C$. Thus, geometric analogues of
unramified automorphic representations should be some geometric
objects which ``live'' on a moduli space of $G$-bundles.

Unfortunately, for a non-abelian group $G$ there is no algebraic
variety whose set of $\C$-points is the set of isomorphism classes of
$G$-bundles on $X$ (for $G=GL_1$ we can take the Picard
variety). However, there is an {\em algebraic moduli stack} denoted by
$\Bun_G$.  It is not an algebraic variety, but it looks locally like
the quotient of an algebraic variety by the action of an algebraic
group (these actions are not free, and therefore the quotient is no
longer an algebraic variety). It turns out that this is good enough
for our purposes.

So which geometric objects on $\Bun_G$ will replace unramified
automorphic representations? Here we need to recall that the function
on the double quotient $G(F)\bs G(\AD_F)/G(\OO_F)$ attached to an
unramified automorphic representation has a special property: it is an
eigenfunction of the so-called Hecke operators. Those are cousins of
the classical Hecke operators one studies in the theory of modular
forms (which is in the left column of Weil's big picture). The
geometric objects we are looking for will be certain sheaves on
$\Bun_G$ satisfying an analogue of the Hecke property. We will call
them {\em Hecke eigensheaves}.

More precisely, these sheaves are ${\mc D}$-{\em modules} on
$\Bun_G$. Recall (see, e.g., \cite{KS,GM}) that a ${\mc D}$-module on
a smooth algebraic variety $Z$ is a sheaf of modules over the sheaf
${\mc D}_Z$ of differential operators on $Z$. An example of a ${\mc
D}$-module is the sheaf of sections of a flat vector bundle on
$Z$. The sheaf of functions on $Z$ acts on sections by multiplication,
so it is an $\OO_Z$-module. But the flat connection also allows us to
act on sections by vector fields on $Z$. This gives rise to an action
of the sheaf ${\mc D}_Z$, because it is generated by vector fields and
functions. Thus, we obtain the structure of a ${\mc D}$-module.

In our case, $\Bun_G$ is not a variety, but an algebraic stack, but
the (derived) category of ${\mc D}$-modules on it has been defined in
\cite{BD}. On this category act the so-called {\em Hecke
functors}. These are labeled by pairs $(x,V)$, where $x \in X$ and $V$
is a finite-dimensional representation of the dual group $\LG$, and
are defined using certain modifications of $G$-bundles.

Instead of giving a general definition (which may be found in
\cite{BD} or \cite{F:houches}) we will consider two examples. First,
consider the abelian case when $G = GL_1$ (thus, we have $G(\C) =
\C^\times$). In this case $\Bun_G$ may be replaced by the Picard
variety $\on{Pic}$ which parametrizes line bundles on $X$. Given a
point $x \in X$, consider the map $h_x: \on{Pic} \to \on{Pic}$ sending
a line bundle ${\mc L}$ to ${\mc L}(x)$ (the line bundle whose
sections are sections of ${\mc L}$ which are allowed to have a pole of
order $1$ at $x$). By definition, the Hecke functor $H_{1,x}$
corresponding to $x$ and $1 \in \Z$ (which we identify with the set of
one-dimensional representations of $\LG = GL_1$), is given by the
formula
$$
H_{1,x}({\mc F}) = h_x^*({\mc F}).
$$

Next, consider the case of $G=GL_n$ and $V=V_{\check\omega_1}$, the
defining $n$-dimensional representation of $\LG = GL_n$. In this case
$\Bun_{GL_n}$ is the moduli stack $\Bun_n$ of rank $n$ bundles on
$X$. There is an obvious analogue of the map $h_x$, sending a rank $n$
bundle ${\mc M}$ to ${\mc M}(x)$. But then the degree of the bundle
jumps by $n$. It is possible to increase it by $1$, but we need
to choose a line $\ell$ in the fiber of ${\mc M}$ at $x$. We then
define a new rank $n$ bundle ${\mc M}'$ by saying that its sections
are the sections of ${\mc M}$ having a pole of order $1$ at $x$, but
the polar part has to belong to $\ell$. Then $\on{deg} {\mc M}' =
\on{deg} {\mc M} + 1$. However, we now have a ${\mathbb P}^{n-1}$
worth of modifications of ${\mc M}$ corresponding to different choices
of the line $\ell$. The Hecke functor $H_{V_{\check\omega_1,x}}$ is
obtained by ``integrating'' over all of them.

More precisely, let ${\mc H}ecke_{\check\omega_1,x}$ be the moduli
stack of pairs $(\M,\M')$ as above. It defines a correspondence over
$\on{Bun}_n \times \on{Bun}_n$:
\begin{equation}    \label{Hecke cor}
\begin{array}{ccccc}
& & {\mc Hecke}_{\check\omega_1,x} & & \\
& \stackrel{\hl_x}\swarrow & & \stackrel{\hr_x}\searrow & \\
\Bun_n & & & & \Bun_n
\end{array}
\end{equation}
By definition,
\begin{equation}    \label{formula H1}
H_{\check\omega_1,x}({\mc F}) = \hr_{x*} \; \hl_x{}^*({\mc
  F}).
\end{equation}

For irreducible representations $V_{\cla}$ of $\LG$ with general
dominant integral highest weights $\cla$ there is an analogous
correspondence in which the role of the projective space ${\mathbb
P}^{n-1}$ is played by the Schubert variety in the affine Grassmannian
of $G$ corresponding to $\cla$ (see \cite{BD,MV}, and \cite{F:houches}
for a brief outline).

Allowing the point $x$ vary, we obtain a correspondence between
$\Bun_G$ and $X \times \Bun_G$ and Hecke functors acting from the
category of ${\mc D}$-modules on $\Bun_G$ to the (derived) category of
${\mc D}$-modules on $X \times \Bun_G$, which we denote by $H_V, V \in
\Rep \LG$.

Now let ${\mc E} = (E,\nabla)$ be a flat $\LG$-bundle on $X$. A ${\mc
D}$-module ${\mc F}$ on $\Bun_G$ is called a {\em Hecke eigensheaf}
with respect to ${\mc E}$ (or with ``eigenvalue'' ${\mc E}$) is we
have a collection of isomorphisms
\begin{equation}    \label{hecke}
H_V({\mc F}) \simeq V_{\mc E} \boxtimes {\mc F},
\end{equation}
compatible with the tensor product structures. Here 
$$
V_{\mc E} = {\mc E} \underset{\LG}\times V
$$
is the flat vector bundle on $X$ associated to ${\mc E}$ and $V$,
viewed as a ${\mc D}$-module. Thus, in particular, we have a
collection of isomorphisms
$$
H_{V,x}({\mc F}) \simeq V \otimes {\mc F}, \qquad x \in X.
$$
When we vary the point $x$, the ``eigenvalues'', which are all
isomorphic to the vector space underlying $V$, combine into the flat
vector bundle $V_{\mc E}$ on $X$.

The {\em geometric Langlands conjecture} may be stated as follows: for
any flat $\LG$-bundle ${\mc E}$ there exists a non-zero ${\mc
D}$-module ${\mc F}_{\mc E}$ on $\Bun_G$ with eigenvalue ${\mc E}$.

Moreover, if ${\mc E}$ is irreducible, this ${\mc D}$-module is
supposed to be irreducible (when restricted to each connected
component of $\Bun_G$) and unique up to an isomorphism (it should also
be holonomic and have regular singularities). But if ${\mc E}$ is not
irreducible, we might have a non-trivial (derived) category of Hecke
eigensheaves, and the situation becomes more subtle.

Thus, at least for irreducible ${\mc E}$, we expect the following
picture:

\begin{equation}    \label{LC}
\boxed{\begin{matrix} \text{flat} \\
    \LG\text{-bundles on } X \end{matrix}} \quad
    \longrightarrow \quad \boxed{\begin{matrix}
    \text{Hecke eigensheaves} \\ \text{on } \Bun_G \end{matrix}}
\end{equation}

$$
{\mc E} \quad \longrightarrow \quad {\mc F}_{\mc E}.
$$

\medskip

The geometric Langlands correspondence has been constructed in many
cases. For $G=GL_n$ the Hecke eigensheaves corresponding to
irreducible ${\mc E}$ have been constructed in \cite{FGV,Ga}, building
on the work of P. Deligne for $n=1$ (explained in \cite{Laumon:duke}
and \cite{F:houches}), V. Drinfeld \cite{Dr1} for $n=2$, and G. Laumon
\cite{Laumon:duke} (this construction works for curves defined both
over $\Fq$ or $\C$).

For all simple algebraic groups $G$ the Hecke eigensheaves have been
constructed in a different way (for curves over $\C$) by A. Beilinson
and V. Drinfeld \cite{BD} in the case when ${\mc E}$ has an additional
structure of an {\em oper} (this means that ${\mc E}$ belongs to a
certain half-dimensional locus in $\Loc_{\LG}$). It is interesting
that this construction is also closely related to quantum field
theory, but in a seemingly different way. Namely, it uses methods of
2D Conformal Field Theory and representation theory of affine
Kac--Moody algebras of critical level. For more on this, see Part III
of \cite{F:houches}.

\section{CATEGORICAL VERSION}    \label{cat ver}

Looking at the correspondence \eqref{LC}, we notice that there is an
essential asymmetry between the two sides. On the left we have flat
$\LG$-bundles, which are points of a moduli stack $\Loc_{\LG}$ of
flat $\LG$-bundles (or local systems) on $X$. But on the right we
have Hecke eigensheaves, which are objects of a category; namely, the
category of ${\mc D}$-modules on $\Bun_G$. Beilinson and Drinfeld have
suggested a natural way to formulate it in a more symmetrical way.

The idea is to replace a point ${\mc E} \in \Loc_{\LG}$ by an object
of another category; namely, the skyscraper sheaf ${\mc O}_{\mc E}$ at
${\mc E}$ viewed as an object of the category of coherent ${\mc
O}$-modules on $\Loc_{\LG}$. A much stronger, categorical, version of
the geometric Langlands correspondence is then a conjectural
equivalence of derived categories\footnote{\label{foot} It is expected
(see \cite{FW}, Sect. 10) that there is in fact a $\Z_2$-gerbe of such
equivalences. This gerbe is trivial, but not canonically
trivialized. One gets a particular trivialization of this gerbe, and
hence a particular equivalence, for each choice of the square root of
the canonical line bundle $K_X$ on $X$.}

\begin{equation}    \label{na fm}
\boxed{\begin{matrix} \text{derived category of} \\
    \OO\text{-modules on } \on{Loc}_{\LG} \end{matrix}} \quad
    \longleftrightarrow \quad \boxed{\begin{matrix}
    \text{derived category of} \\ \D\text{-modules on } \Bun_G
    \end{matrix}}
\end{equation}

\bigskip

This equivalence should send the skyscraper sheaf ${\mc O}_{\mc E}$ on
$\on{Loc}_{\LG}$ supported at ${\mc E}$ to the Hecke eigensheaf ${\mc
F}_E$. If this were true, it would mean that Hecke eigensheaves
provide a good ``basis'' in the category of $\D$-modules on $\Bun_G$,
so we would obtain a kind of spectral decomposition of the derived
category of $\D$-modules on $\Bun_G$, like in the Fourier
transform. (Recall that under the Fourier transform on the real line
the delta-functions $\delta_x$, analogues of ${\mc O}_{\mc E}$, go to
the exponential functions $e^{itx}$, analogues of ${\mc F}_{\mc E}$.)

This equivalence has been proved by G. Laumon \cite{Laumon:F} and
M. Rothstein \cite{Rothstein} in the abelian case, when $G=GL_1$ (or a
more general torus). They showed that in this case this is nothing but
a version of the Fourier--Mukai transform. Thus, the categorical
Langlands correspondence may be viewed as a kind of non-abelian
Fourier--Mukai transform (see \cite{F:houches}, Section 4.4).

Unfortunately, a precise formulation of such a correspondence, even as
a conjecture, is not so clear because of various subtleties
involved. One difficulty is the structure of $\on{Loc}_{\LG}$. Unlike
the case of $\LG = GL_1$, when all flat bundles have the same groups
of automorphisms (namely, $GL_1$) and $\Loc_{GL_1}$ is smooth, for a
general group $\LG$ the groups of automorphisms are different for
different flat bundles, and so $\on{Loc}_{\LG}$ is a complicated
stack. For example, if $\LG$ is a simple Lie group of adjoint type,
then a generic flat $\LG$-bundle has no automorphisms, while the group
of automorphisms of the trivial flat bundle is isomorphic to
$\LG$. In addition, unlike $\Bun_G$, the stack $\Loc_{\LG}$ has
singularities. All of this has to be reflected on the other side of
the correspondence, in ways that have not yet been fully understood.

Nevertheless, the diagram \eqref{na fm} gives us a valuable guiding
principle to the geometric Langlands correspondence. In particular, it
gives us a natural explanation as to why the skyscraper sheaves on
$\Loc_{\LG}$ should correspond to Hecke eigensheaves.

The point is that on the category of ${\mc O}$-modules on $\Loc_{\LG}$
we also have a collection of functors $W_{V}$, parametrized by the
same data as the Hecke functors $H_{V}$. Following physics
terminology, we will call them {\em Wilson functors}. These functors
act from the category of $\OO$-modules on $\on{Loc}_{\LG}$ to the
category of sheaves on $X \times \on{Loc}_{\LG}$, which are ${\mc
D}$-modules along $X$ and $\OO$-modules along $\on{Loc}_{\LG}$.

To define them, observe that we have a tautological $\LG$-bundle
${\mc T}$ on $X \times \on{Loc}_{\LG}$, whose restriction to $X
\times {\mc E}$, where ${\mc E} = (E,\nabla)$, is $E$. Moreover,
$\nabla$ gives us a partial connection on ${\mc T}$ along $X$. For a
representation $V$ of $\LG$, let $V_{\mc T}$ be the associated vector
bundle on $X \times \on{Loc}_{\LG}$, with a connection along $X$.

Let $p: X \times \on{Loc}_{\LG} \to \on{Loc}_{\LG}$ be the projection
onto the second factor. By definition,
\begin{equation}    \label{wilson}
W_V({\mc F}) = V_{\mc T} \otimes p^*({\mc F}).
\end{equation}
(note that by construction $V_{\mc T}$ carries a connection along $X$
and so the right hand side really is a ${\mc D}$-module along $X$).

Now, the conjectural equivalence \eqref{na fm} should be compatible
with the Wilson/Hecke functors in the sense that
\begin{equation}    \label{WH}
C(W_V({\mc F})) \simeq H_V(C({\mc F})), \qquad V \in \on{Rep} {}\LG,
\end{equation}
where $C$ denotes this equivalence (from left to right).

In particular, observe that the skyscraper sheaf $\OO_{\mc E}$ at
${\mc E} \in \on{Loc}_{\LG}$ is obviously an eigensheaf of the Wilson
functors:
$$
W_V(\OO_{\mc E}) = V_{\mc E} \boxtimes \OO_{\mc E}.
$$
Indeed, tensoring a skyscraper sheaf with a vector bundle is
the same as tensoring it with the fiber of this vector bundle at the
point of support of this skyscraper sheaf. Therefore \eqref{WH}
implies that ${\mc F}_{\mc E} = C(\OO_{\mc E})$ must satisfy the
Hecke property \eqref{hecke}. In other words, ${\mc F}_{\mc E}$ should
be a Hecke eigensheaf on $\Bun_G$ with eigenvalue ${\mc E}$. Thus, we
obtain a natural explanation of the Hecke property of ${\mc F}_{\mc
  E}$: it follows from the compatibility of the categorical Langlands
correspondence \eqref{na fm} with the Wilson/Hecke functors.

Let us summarize: the conjectural equivalence \eqref{na fm} gives us a
natural and convenient framework for the geometric Langlands
correspondence. It is this equivalence that Kapustin and Witten have
related to the $S$-duality of 4D super--Yang--Mills.

\section{ENTER PHYSICS}    \label{PHY}

We will now add a fourth column to Weil's big picture, which we will
call ``Quantum Physics'':

$$
\begin{matrix}
\text{{\em Number Theory}} \quad & \quad \text{{\em Curves over}
  $\Fq$} \quad & \quad \text{{\em Riemann Surfaces}} \quad & \quad
  \text{{\em Quantum Physics}}
\end{matrix}
$$

\medskip

In the context of the Langlands Program, the last column means
$S$-duality and Mirror Symmetry of certain 4D and 2D quantum field
theories, which we will now briefly describe following \cite{KW}.

We start with the pure 4D Yang--Mills (or gauge) theory on a
Riemannian four-manifold $M_4$. Let $G_c$ be a compact connected
simple Lie group. The classical (Euclidean) action is a functional on
the space of connections on arbitrary principal $G_c$-bundles ${\mc
P}$ on $M_4$ given by the formula
$$
I = \frac{1}{4g^2} \int_{M_4} \on{Tr} F_A \wedge \star F_A + \frac{i
  \theta}{8\pi^2} \int_{M_4} \on{Tr} F_A \wedge F_A.
$$

\noindent Here $F_A$ is the curvature of the connection $A$ (a
$\g$-valued two-form on $M_4$), $\star\,$ is the Hodge star operator,
and $\on{Tr}$ is the invariant bilinear form on the Lie algebra $\g$
normalized in such a way that the second term is equal to $i\theta k$,
where $k$ could be an arbitrary integer, if $G_c$ is
simply-connected. The second term is equal to $i\theta$ times the
second Chern class $c_2({\mc P})$ of the bundle ${\mc P}$ and hence is
topological. Correlation functions are given by path integrals of the
form $\int e^{-I}$ over the space of connections modulo gauge
transformations. Hence they may be written as Fourier series in
$e^{i\theta}$ (or its root if $G_c$ is not simply-connected) such that
the coefficient in front of $e^{i\theta n}$ is the sum of
contributions from bundles ${\mc P}$ with $c_2({\mc P}) = -n$.

It is customary to combine the two parameters, $g$ and $\theta$, into
one complex coupling constant
$$
\tau = \frac{\theta}{2\pi} + \frac{4\pi i}{g^2}.
$$

Next, we consider $N=4$ supersymmetric extension of this model. This
means that we add fermionic and bosonic fields in such a way that the
action of the Lorentz group (we will work in Euclidean signature,
where this group becomes $SO(4)$) is extended to an action of an
appropriate supergroup (see \cite{KW} or the books \cite{Phys} for
background on supersymmetric quantum field theory).

The $S$-duality of this theory is the statement that the theory with
gauge group $G_c$ and complex coupling constant $\tau$ is equivalent
to the theory with the Langlands dual gauge group $\LG_c$ and
coupling constant $^L \tau = - 1/n_{\g}\tau$:
\begin{equation}    \label{S-duality1}
(G_c,\tau) \longleftrightarrow ({}\LG_c,-1/n_{\g} \tau),
\end{equation}
where $n_{\g}$ is the lacing number of the Lie algebra $\g$ (equal to
$1$ for simply-laced Lie algebras, $2$ for $B_n,C_n$ and $F_4$, and
$3$ for $G_2$). This is an extension of the duality \eqref{S-duality}
of \cite{MO} discussed in the introduction to non-zero values of
$\theta$ (with $g$ normalized in a slightly different way). In
addition, for simply-laced $G_c$ the path integral is a Fourier series
in $e^{i\theta}$, so $\theta$ may be shifted by an integer multiple of
$2\pi$ without changing the path integral. Thus, we also have the
equivalence
$$
(G_c,\tau) \longleftrightarrow (G_c,\tau+1).
$$
For general non-simply connected groups we have instead a symmetry
$\tau \mapsto \tau + n\Z$, where $n$ is a certain
integer. Thus, we obtain an action of a subgroup of $SL_2(\Z)$ on the
super--Yang--Mills theories with gauge groups $G_c$ and
$\LG_c$.\footnote{In general, it is a proper subgroup of $SL_2(\Z)$
for two reasons: first, we have the coefficient $n_{\g}$ in formula
\eqref{S-duality1} for non-simply laced $G_c$, and second, the dual of
a simply-connected Lie group is not simply-connected in general, in
which case the transformation $\tau \to \tau+1$ is not a symmetry.}
As we discussed in the Introduction, this is a striking statement
because it relates a theory at strong coupling to a theory at weak
coupling.

We want to focus next on a ``topological sector'' of this theory. This
means that we pick an element $Q$ in the Lie superalgebra ${\mathfrak
s}$ of the super-Lorentz group (the supergroup extension of $SO(4)$)
such that $Q^2=0$, and such that the stress tensor (which is a field
responsible for variation of the metric on $M_4$) is equal to the
commutator of $Q$ and another field. Let us restrict ourselves to
those objects (fields, boundary conditions, etc.) in the theory which
commute with this $Q$. This is a particular (and relatively small)
sector of the full quantum field theory, in which all quantities (such
as correlation functions) are topological, that is,
metric-independent. This sector is what is usually referred to as
Topological Field Theory (TFT).

There is a problem, however. For this $Q \in {\mathfrak s}$ to be
well-defined on an arbitrary manifold $M_4$, it has to be invariant
under the action of the Lorentz group $SO(4)$ -- more precisely, its
double cover $Spin(4)$. Unfortunately, there are no such elements in
our Lie superalgebra ${\mathfrak s}$. In order to obtain such an
element, one uses a trick, called {\em twisting} (see, e.g.,
\cite{Witten:atiyah}). Our theory has an additional group of
automorphisms commuting with the action of $Spin(4)$, called
$R$-symmetry; namely, the group $Spin(6)$. We can use it to modify the
action of $Spin(4)$ on the fields of the theory and on the Lie
superalgebra ${\mathfrak s}$ as follows: define a new action of
$Spin(4)$ equal to the old action together with the action coming from
a homomorphism $Spin(4) \to Spin(6)$ and the action of $Spin(6)$ by
$R$-symmetry. One might then be able to find a differential $Q \in
{\mathfrak s}$ invariant under this new action of $Spin(4)$.

There are essentially three different choices for doing this, as
explained in \cite{VW}. The first two are similar to the twists used
in Witten's construction of a topological field theory that yields
Donaldson invariants of four-manifolds (which is a topological twist
of an $N=2$ supersymmetric Yang--Mills theory) \cite{Witten:don}. It
is the third twist, studied in detail in \cite{KW}, that is relevant
to the geometric Langlands. For this twist there are actually two
linearly independent (and anti-commuting with each other) operators,
$Q_l$ and $Q_r$, which square to $0$. We can therefore use any linear
combination
$$
Q = u Q_l + v Q_r
$$
as the differential defining the topological field theory (for each of
them the stress tensor will be a $Q$-commutator, so we will indeed
obtain a topological field theory). We obtain the same theory if we
rescale both $u$ and $v$ by the same number. Hence we obtain a family
of topological field theories parametrized by ${\mathbb P}^1$.

Let $t=v/u$ be a coordinate on this ${\mathbb P}^1$. We will refer to
it as the ``twisting parameter.''  The $S$-duality \eqref{S-duality1}
should be accompanied by the change of the twisting parameter
according to the rule
\begin{equation}    \label{t}
t \mapsto \frac{\tau}{|\tau|} t.
\end{equation}

How can we test $S$-duality of these topological theories? Vafa and
Witten have earlier tested the $S$-duality of a different (Donaldson
type) topological twisting of $N=4$ super--Yang--Mills by showing that
the partition functions of these theories (depending on $\tau$) are
modular forms \cite{VW}. This proves the invariance of the partition
functions under the action of a subgroup of $SL_2(\Z)$ on $\tau$. What
turns out to be relevant to the geometric Langlands Program is the
study of boundary conditions in these topological field theories.

Kapustin and Witten assume that the four-manifold $M_4$ has the form
$$
M_4 = \Sigma \times X,
$$
where $X$ is a closed Riemann surface (this will be the algebraic
curve of the geometric Langlands) and $\Sigma$ is a Riemann surface
with a boundary (which we may simply take to be a half-plane). They
study the limit of the topological gauge theory on this manifold when
$X$ becomes very small (this is called ``compactification of the
theory on $X$''). In this limit the theory is described by an
effective {\em two-dimensional} topological field theory on
$\Sigma$. In earlier works \cite{BJSV,HMS} the latter theory was
identified with the (twisted) topological sigma model on $\Sigma$ with
the target manifold ${\mc M}_H(G)$, the {\em Hitchin moduli space} of
Higgs $G$-bundles on $X$. Moreover, the $S$-duality of the
supersymmetric gauge theories on $\Sigma \times X$ (for particular
values of $\tau$ and $t$) becomes Mirror Symmetry between the
topological sigma models with the targets ${\mc M}_H(G)$ and ${\mc
M}_H({}\LG)$.

Next, we look at the boundary conditions in the gauge theories, which
give rise to {\em branes} in these sigma models. $S$-duality yields an
equivalence of the categories of branes for ${\mc M}_H(G)$ and ${\mc
M}_H({}\LG)$ (also known, after M. Kontsevich, as Homological Mirror
Symmetry). Kapustin and Witten have related this equivalence to the
categorical geometric Langlands correspondence \eqref{na fm}. Thus,
they establish a link between $S$-duality and geometric Langlands
duality. We describe this in more detail in the next section.

\section{MIRROR SYMMETRY OF HITCHIN MODULI SPACES}    \label{MIRROR}

In \cite{Hit} N. Hitchin introduced a remarkable hyper-K\"ahler
manifold ${\mc M}_H(G)$ for each smooth projective complex algebraic
curve $X$ and reductive Lie group $G$. It is easiest to describe it in
its complex structure $I$, in which it is the moduli space of
semi-stable {\em Higgs bundles} on $X$. Recall that a Higgs $G$-bundle
on $X$ is a pair $(E,\phi)$, where $E$ is a (algebraic) $G$-bundle on
$X$ and $\phi$ is a Higgs field on it, that is,
$$
\phi \in H^0(X,\g_E \otimes K_X),
$$
where $\g_E = E \underset{G}\times \g$ is the adjoint vector bundle.

In the complex structure $J$, however, ${\mc M}_H(G)$ is described as
the moduli space of semi-stable {\em flat bundles}, that is, pairs
$(E,\nabla)$, where $E$ is again (algebraic) $G$-bundle on $X$ and
$\nabla$ is (algebraic) connection on $E$. To distinguish between it
as a complex algebraic variety from the moduli space of Higgs bundles
we will denote it by ${\mc Y}(G)$. The two are isomorphic as real
manifolds (this is the statement of non-abelian Hodge theory
\cite{Hit,Corlette,Simpson}), but not as complex (or algebraic)
manifolds.

There are two types of twisted supersymmetric two-dimensional sigma
models with K\"ahler target manifolds: $A$-model and $B$-model (see
\cite{Witten:mirror}). The former depends on the symplectic structure
on the target manifold and the latter depends on the complex
structure.

Kapustin and Witten start with two topological twisted
super--Yang--Mills theories on $\Sigma \times X$. One has gauge group
$G_c$, twisting parameter $t=1$, and $\theta=0$. The other, $S$-dual
theory, has gauge group $\LG_c$, the twisting parameter $^L t = i$,
and $^L \theta = 0$ (neither of these topological theories depends on
$g$ \cite{KW}).\footnote{As explained in \cite{KW,K}, for some other
values of parameters one obtains the so-called quantum geometric
Langlands correspondence (see \cite{F:houches}, Section 6.3).} They
show that after compactification on $X$ the first theory becomes
the $A$-model with the target manifold ${\mc M}_H(G)$ and the
symplectic structure $\omega_K$, which is the K\"ahler form for the
complex structure $K$ on ${\mc M}_H(G)$. This symplectic structure has
a nice geometric description. Note that the Higgs field $\phi$ is an
element of $H^0(X,\g_E \otimes K_X)$, which is isomorphic to the
cotangent space to $E$, viewed as a point of $\Bun_G$, the moduli
stack of $G$-bundles on $X$. Thus, ${\mc M}_H(G)$ is almost the
cotangent bundle to $\Bun_G$; ``almost'' because we impose the
semi-stability condition on the Higgs bundle. The symplectic form
$\omega_K$ comes from the standard symplectic form on the cotangent
bundle (which is the imaginary part of the holomorphic symplectic
form).

The second gauge theory becomes, after compactification on $X$, the
$B$-model with the target manifold ${\mc Y}({}\LG)$; that is, ${\mc
M}_H({}\LG)$ with respect to the complex structure $J$.

After dimensional reduction from 4D to 2D, the $S$-duality of
super--Yang--Mills theories becomes Mirror Symmetry between the
$A$-model with the target manifold ${\mc M}_H(G)$ (and symplectic
structure $\omega_K$) and the $B$-model with the target manifold ${\mc
Y}({}\LG)$ (and complex structure $J$).

\medskip

\noindent {\em Remark.} As explained in \cite{KW}, there is also
Mirror Symmetry between $A$- and $B$-models with respect to other
symplectic and complex structures. For instance, there is Mirror
Symmetry, studied in \cite{DP,Hi3,Ar}, between the $B$-models on ${\mc
M}_H(G)$ and ${\mc M}_H(\LG)$ with respect to the complex structures
$I$ on both of them. In what follows we will not discuss these
additional dualities.

\subsection{Dual Hitchin fibrations} In order to understand better
this Mirror Symmetry, we recall the construction of the {\em Hitchin
map}. For any Higgs bundle $(E,\phi)$ and an invariant polynomial $P$
of degree $d$ on the Lie algebra $\g$, we can evaluate $P$ on $\phi$
and obtain a well-defined section $P(\phi)$ of $K_X^{\otimes d}$. The
algebra $(\on{Fun}(\g))^G$ of invariant polynomial functions on $\g$
is a graded free polynomial algebra with $\ell=\on{rank}(\g)$
generators of degrees $d_i, i=1,\ldots,\ell$, where the $d_i$ are the
exponents of $\g$ plus $1$. Let us choose a set of generators $P_i,
i=1,\ldots,\ell$. Then we construct the Hitchin map \cite{Hit,Hi2}
\begin{align*}
p: {\mc M}_H(G) \to {\mb B} &= \bigoplus_{i=1}^\ell
H^0(X,K_X^{\otimes d_i}), \\
(E,\phi) &\mapsto (P_1(\phi),\ldots,P_\ell(\phi)).
\end{align*}
This is slightly non-canonical, because there is no canonical choice
of generators $P_i$ in general. More canonically, we have a map to
$$
{\mb B} := H^0(X,(\g/\!/G)_{K_X}), \qquad (\g/\!/G)_{K_X} := K_X^\times
\underset{\C^\times}\times \g/\!/G,
$$
where $K_X^\times$ denotes the $\C^\times$-bundle associated to $K_X$,
and $\g/\!/G := \on{Spec}((\on{Fun}(\g))^G)$ is the graded vector space
on which the $P_i$ are coordinate functions (the $\C^\times$-action on
it comes from the grading).

By Chevalley's theorem, $\g/\!/G = \h/\!/W :=
\on{Spec}((\on{Fun}(\h))^W)$, where $\h$ is a Cartan subalgebra. By
definition, $^L\h = \h^*$. Hence $^L\g/\!/{}\LG = {}^L\h/\!/W =
\h^*/\!/W$. Choosing a non-degenerate invariant bilinear form
$\kappa_0$ on $\g$, we identify $\h$ and $\h^*$, and hence ${\mb B}$
and $^L {\mb B}$. Any other invariant bilinear form is proportional to
$\kappa_0$. Hence replacing $\kappa_0$ by another non-zero bilinear
form would correspond to a $\C^\times$-action on the base. But this
action can be lifted to a $\C^\times$-action on the total space ${\mc
M}_H(G)$ (rescaling the Higgs field $\phi$). Hence the ambiguity in
the choice of $\kappa_0$ is not essential; it may be absorbed into an
automorphism of one of the two Hitchin moduli spaces.

As the result, we obtain two fibrations over the same base ${\mb B}$:

\begin{equation}    \label{mirror}
\begin{array}{ccccc}
{\mc M}_{H}(\LG) & \; & \; & \; & {\mc M}_H(G) \\
\; & \searrow & \; & \swarrow & \; \\
\; & \; & {\mb B} & \; & \;
\end{array}
\end{equation}

\medskip

For generic $b \in {\mb B}$ (the connected components of) the fibers
$^L {\mb F}_b$ and ${\mb F}_b$ of these Hitchin fibrations are smooth
tori, which are in fact isomorphic to abelian varieties in the complex
structure $I$. For instance, for $G=SL_n$, ${\mb F}_b$ is the
generalized Prym variety of the spectral curve associated to $b$,
which is a smooth degree $n$ cover of $X$ if $b$ is generic.

Moreover, the tori $^L {\mb F}_b$ and ${\mb F}_b$ (again, for generic
$b$) are {\em dual} to each other. This can be expressed in the
following way which will be convenient for our purposes: there is a
bijection between points of $^L {\mb F}_b$ and flat unitary line
bundles on ${\mb F}_b$.\footnote{\label{foot1} This bijection depends
on the choice of base points in the fibers. Using the Hitchin section,
one obtains such base points, but for this one may have to choose a
square root of the canonical line bundle $K_X$. This is closely
related the $\Z_2$-gerbe ambiguity in the equivalence \eqref{na fm}
discussed in the footnote on page \pageref{foot}.}  The duality of
tori is the simplest (abelian) example of Mirror Symmetry, also known
as $T$-{\em duality}. Thus, we expect that the Mirror Symmetry of the
Hitchin moduli spaces ${\mc M}_{H}(\LG)$ and ${\mc M}_H(G)$ is
realized via fiberwise $T$-duality (for generic fibers of the dual
Hitchin fibrations). This is an example of a general proposal of
Strominger, Yau and Zaslow \cite{SYZ} that Mirror Symmetry of two
Calabi--Yau manifolds $X$ and $Y$ should be realized as $T$-duality of
generic fibers in special Lagrangian fibrations of $X$ and $Y$ (what
happens for the singular fibers is {\em a priori} far less clear; but
see \secref{endoscopy}). It is in this sense that ``$S$-duality
reduces to $T$-duality'' \cite{HMS}.

We are interested in the study of the $B$-model with the target ${\mc
M}_H({}\LG)$, with respect to the complex structure $J$ (that is, the
moduli space ${\mc Y}({}\LG)$ of flat bundles), and the $A$-model with
the target ${\mc M}_H(G)$, with respect to the symplectic structure
$\omega_K$. These two topological field theories are expected to be
equivalent to each other. Therefore anything we can say about one of
them should have a counterpart in the other. For instance, their
cohomologies may be interpreted as the spaces of vacua in these field
theories, and hence they should be isomorphic. This has indeed been
verified by Hausel and Thaddeus \cite{HT} in the case when $G=SL_n,
{}\LG = PGL_n$ (since the Hitchin moduli spaces are non-compact,
special care has to be taken to properly define these cohomologies,
see \cite{HT}).

To make contact with the geometric Langlands correspondence, Kapustin
and Witten study in \cite{KW} the {\em categories of branes} in these
two topological field theories.

\subsection{Categories of branes}

Branes in two-dimensional sigma models are certain generalizations of
boundary conditions. When writing path integral for maps $\Phi: \Sigma
\to M$, where $\Sigma$ has a boundary, we need to specify boundary
conditions for $\Phi$ on $\pa \Sigma$. We may also ``couple'' the
sigma model to another quantum field theory on $\pa \Sigma$ (that is,
modify the action by a boundary term) which may be interpreted as a
decoration of the boundary condition. In topological field theory
these conditions should preserve the supersymmetry, which leads to
natural restrictions.

A typical example of a boundary condition is specifying that
$\Phi(\partial \Sigma)$ belongs to a submanifold $M' \subset M$. In
the $B$-model the target manifold $M$ is a complex manifold, and in
order to preserve the supersymmetry $M'$ has to be a complex
submanifold. In the $A$-model, $M$ is a symplectic manifold and $M'$
should be Lagrangian. Coupling to field theories on $\pa \Sigma$
allows us to introduce into the picture a holomorphic vector bundle on
$M'$ in the case of $B$-model, and a flat unitary vector bundle on
$M'$ in the case of $A$-model.

More generally, the category of branes in the $B$-model with a complex
target manifold $M$ (called $B$-branes) is the (derived) category of
coherent sheaves on $M$, something that is fairly well understood
mathematically. The category of branes in the $A$-model with a
symplectic target manifold $M$ (called $A$-branes) is less
understood. It is believed to contain what mathematicians call the
Fukaya category, typical objects of which are pairs $(L,\nabla)$,
where $L \subset M$ is a Lagrangian submanifold and $\nabla$ is a flat
unitary vector bundle on $L$. However (and this turns out to be
crucial for applications to the geometric Langlands), it also contains
more general objects, such as coisotropic submanifolds of $M$ equipped
with vector bundles with unitary connection.

Under the Mirror Symmetry between the sigma models with the target
manifolds ${\mc Y}({}\LG)$ and ${\mc M}_H(G)$ we therefore expect to
have the following equivalence of (derived) categories of branes
(often referred to, after Kontsevich, as Homological Mirror Symmetry):

\begin{equation}    \label{branes}
\boxed{\text{$B$-branes on } {\mc Y}(\LG)} \quad
    \longleftrightarrow \quad \boxed{\text{$A$-branes on } {\mc
    M}_H(G)}
\end{equation}

\bigskip

It is this equivalence that Kapustin and Witten have related to the
categorical Langlands correspondence \eqref{na fm}. The category on
the left in \eqref{branes} is the (derived) category of coherent
sheaves on ${\mc Y}(\LG)$, which is the moduli space of semi-stable
flat $\LG$-bundles on $X$. It is closely related to the (derived)
category of coherent sheaves (or, equivalently, $\OO$-modules) on
$\Loc_{\LG}$, which appears on the left of \eqref{na fm}. The
difference is that, first of all, $\Loc_{\LG}$ is the moduli {\em
stack} of flat $\LG$-bundles on $X$, whereas ${\mc Y}(\LG)$ is the
moduli {\em space} of semi-stable ones. Second, from the physics
perspective it is more natural to consider coherent sheaves on ${\mc
Y}(\LG)$ with respect to its complex analytic rather than algebraic
structure, whereas in \eqref{na fm} we consider algebraic
$\OO$-modules on $\Loc_{\LG}$. These differences aside, these two
categories are very similar to each other. They certainly share many
objects, such as skyscraper sheaves supported at points
corresponding to stable flat $\LG$-bundles which we will discuss
momentarily.

\subsection{Triangle of equivalences} The categories on the right in
\eqref{na fm} and \eqref{branes} appear at first glance to be quite
different. But Kapustin and Witten have suggested that they should be
equivalent to each other as well. Thus, we obtain the following
triangle of derived categories:

\begin{equation}    \label{triangle}
\xymatrix{& \boxed{\text{$A$-branes on } {\mc M}_H(G)} \ar[dd] \\
\boxed{\text{$B$-branes on } {\mc Y}(\LG)} \ar[ur] \ar[dr] \\
& \boxed{{\mc D}\text{-modules on } \Bun_G}
}
\end{equation}

\medskip

The upper arrow represents Homological Mirror Symmetry \eqref{branes}
whereas the lower arrow represents the categorical Langlands
correspondence \eqref{na fm}.

According to \cite{KW}, Section 11, the vertical arrow is another
equivalence that has nothing to do with either Mirror Symmetry or
geometric Langlands. It should be a general statement linking the
(derived) category of ${\mc D}$-modules on a variety $M$ and the
(derived) category of $A$-branes on its cotangent bundle $T^*M$ (recall
that ${\mc M}_H(G)$ is almost equal to $T^* \Bun_G$). Kapustin and
Witten have proposed the following functor from the category of
$A$-branes on $T^*M$ with respect to the symplectic structure ${\rm
Im}\,\Omega$ (where $\Omega$ is the holomorphic symplectic form on
$T^*M$) to the category of ${\mc D}$-modules on $M$:
\begin{equation}    \label{Hom}
{\mc A} \mapsto \on{Hom}({\mc A}_{cc},{\mc A}),
\end{equation}
where ${\mc A}_{cc}$ is a ``canonical coisotropic brane'' on $T^*M $.
This is $T^*M$ itself (viewed as a coisotropic submanifold) equipped
with a line bundle with connection satisfying special properties. They
argued on physical grounds that the right hand side of \eqref{Hom} may
be ``sheafified'' along $M$, and moreover that the corresponding sheaf
of rings $\on{Hom}({\mc A}_{cc},{\mc A}_{\on{cc}})$ is nothing but
the sheaf of differential operators on ${\mc M}_H(G)$.\footnote{More
precisely, it is the sheaf of differential operators acting on a
square root of the canonical line bundle on $\Bun_G$, but we will
ignore this subtlety here.} Hence the (sheafified) right hand side of
\eqref{Hom} should be a ${\mc D}$-module. While this argument has not
yet been made mathematically rigorous, it allows one to describe
important characteristics of the ${\mc D}$-module associated to an
$A$-brane, such as its reducibility, the open subset of $M$ where it is
represented by a local system, the rank of this local systems, and
even its monodromy (see Section 4 of \cite{FW}).

An alternative (and mathematically rigorous) approach to establishing
an equivalence between the categories of $A$-branes and ${\mc
D}$-modules has also been proposed by D. Nadler and E. Zaslow
\cite{NZ,Nad}.

Though a lot of work still needs to be done to distill this connection
and reconcile different approaches, this is clearly a very important
and beautiful idea on its own right.

Thus, according to Kapustin and Witten, the (categorical) geometric
Langlands correspondence \eqref{na fm} may be obtained in two
steps. The first step is the Homological Mirror Symmetry
\eqref{branes} of the Hitchin moduli spaces for two dual groups, and
the second step is the above link between the $A$-branes and ${\mc
D}$-modules.

There is actually more structure in the triangle \eqref{triangle}. On
each of these three categories we have an action of certain functors,
and all equivalences between them are supposed to commute with these
functors. We have already described the functors on two of these
categories in \secref{cat ver}: these are the Wilson and Hecke
functors. The functors acting on the categories of $A$-branes, were
introduced in \cite{KW} as the two-dimensional shadows of the 't Hooft
loop operators in 4D super--Yang--Mills theory. Like the Hecke
functors, they are defined using modifications of $G$-bundles, but
only those modifications which preserve the Higgs field. The
$S$-duality of super--Yang--Mills theories is supposed to exchange the
't Hooft operators and the Wilson operators (whose two-dimensional
shadows are the functors described in \secref{cat ver}), and this is
the reason why we expect the equivalence \eqref{branes} to commute
with the action of these functors.

As explained above, the central objects in the geometric Langlands
correspondence are Hecke eigensheaves attached to flat
$\LG$-bundles. Recall that the Hecke eigensheaf ${\mc F}_{\mc E}$ is
the ${\mc D}$-module attached to the skyscraper sheaf $\OO_{\mc E}$
supported at a point ${\mc E}$ of $\Loc_{\LG}$ under the conjectural
equivalence \eqref{na fm}. These ${\mc D}$-modules have very
complicated structure. What can we learn about them from the point of
view of Mirror Symmetry?

Let us assume first that ${\mc E}$ has no automorphisms other than
those coming from the center of $\LG$. Then it is a smooth point of
${\mc Y}(\LG)$. These skyscraper sheaves $\OO_{\mc E}$ are the
simplest examples of $B$-branes on ${\mc Y}(\LG)$ (called
$0$-branes). What is the corresponding $A$-brane on ${\mc M}_H(G)$?

The answer is surprisingly simple. Let $b \in {\mb B}$ be the
projection of ${\mc E}$ to the base of the Hitchin fibration. For
${\mc E}$ satisfying the above conditions the Hitchin fiber $^L {\mb
F}_b$ is a smooth torus (it is actually an abelian variety in the
complex structure $I$, but now we look at it from the point of view of
complex structure $J$, so it is just a smooth torus). It is identified
(possibly, up to a choice of the square root of the canonical line
bundle $K_X$, see the footnote on page \pageref{foot1}) with the
moduli space of flat unitary line bundles on the dual Hitchin fiber
${\mb F}_b$, which happens to be a Lagrangian submanifold of
$\MH(G)$. The Mirror Symmetry sends the $B$-brane $\OO_{\mc E}$ to the
$A$-brane which is the pair $({\mb F}_b,\nabla_{\mc E})$, the
Lagrangian submanifold ${\mb F}_b$ of ${\mc M}_H(G)$, together with the
flat unitary line bundle on it corresponding to ${\mc E}$:
$$
\OO_{\mc E} \mapsto ({\mb F}_b,\nabla_{\mc E}).
$$

Since $\OO_{\mc E}$ is obviously an eigenbrane of the Wilson functors
(as we discussed in \secref{cat ver}), the $A$-brane $({\mb
F}_b,\nabla_{\mc E})$ should be an eigenbrane of the 't Hooft
functors. This may in fact be made into a precise mathematical
conjecture, and Kapustin and Witten have verified it explicitly in
some cases.

Thus, the $A$-branes associated to the simplest $B$-branes turn out to
be very nice and simple. This is in sharp contrast with the structure
of the corresponding ${\mc D}$-modules, which is notoriously
complicated in the non-abelian case. Therefore the formalism of
$A$-branes developed in \cite{KW} has clear advantages. It replaces
${\mc D}$-modules with $A$-branes that are much easier to ``observe
experimentally'' and to analyze explicitly. One can hope to use this
new language in order to gain insights into the structure of the
geometric Langlands correspondence. It has already been used in
\cite{FW} for understanding what happens in the endoscopic case as
explained in the next section.

\subsection{Ramification}    \label{ram}

Up to now we have considered the {\em unramified} case of the
geometric Langlands correspondence, in which the objects on the Galois
side of the correspondence are holomorphic $\LG$-bundles on our curve
$X$ with a holomorphic connection. These flat bundles give rise to
homomorphisms $\pi_1(X) \to \LG$. In the classical Langlands
correspondence one looks at more general homomorphisms $\pi_1(X \bs \{
x_1,\ldots,x_n \}) \to \LG$. Thus, we look at holomorphic
$\LG$-bundles on $X$ with {\em meromorphic} connections which have
poles at finitely many points of $X$. The connections with poles of
order one (regular singularities) correspond to {\em tame}
ramification in the classical Langlands Program. Those with poles of
orders higher than one (irregular singularities) correspond to {\em
wild} ramification.

Mathematically, the ramified geometric Langlands correspondence has
been studied in \cite{FG} and follow-up papers (see \cite{F:book} for
an exposition), using the affine Kac--Moody algebras of critical level
and generalizing the Beilinson--Drinfeld approach \cite{BD} to allow
ramification.

S. Gukov and E. Witten \cite{GW} have explained how to include tame
ramification in the $S$-duality picture. Physicists have a general way
of including into a quantum field theory on a manifold $M$ objects
supported on submanifolds of $M$. An example of this is the {\em
surface operators} in 4D super--Yang--Mills theory, supported on
two-dimensional submanifolds of the four-manifold $M_4$. If we
include such an operator, we obtain a certain modification of the
theory. Let us again take $M_4 = \Sigma \times X$ and take this
submanifold to be of the form $\Sigma \times x, x \in X$. Gukov and
Witten show that for a particular class of surface operators the
dimensional reduction of the resulting theory is the sigma model on
$\Sigma$ with a different target manifold $\M_H(G,x)$, the moduli
space of semi-stable Higgs bundles with regular singularity at $x \in
X$. It is again hyper-K\"ahler, and in the complex structure $J$ it
has a different incarnation as the moduli space of semi-stable bundles
with a connections having regular singularity (see \cite{Simp}).

The moduli space $\M_H(G,x)$ has parameters $(\al,\beta,\gamma)$,
which lie in the (compact) Cartan subalgebra of $\g$ (see
\cite{Simp}). For generic parameters, this moduli space parametrizes
semi-stable triples $(E,\phi,{\mc L})$, where $E$ is a (holomorphic)
$G$-bundle, $\phi$ is a Higgs field which has a pole at $x$ of order
one whose residue belongs to the regular semi-simple conjugacy class
of $\frac{1}{2}(\beta + i \gamma)$, and ${\mc L}$ is a flag in the
fiber of $E$ at $x$ which is preserved by this residue (the remaining
parameter $\alpha$ determines the flag). Various degenerations of
parameters give rise to similar moduli spaces in which the residues of
the Higgs fields could take arbitrary values. The moduli spaces ${\mc
M}_H(G,x)$ and ${\mc M}_H(\LG,x)$ (with matching parameters) are
equipped with a pair of mirror dual Hitchin fibrations, and the Mirror
Symmetry between them is again realized as fiberwise $T$-duality (for
generic fibers which are again smooth dual tori).

The $S$-duality of the super--Yang--Mills theories, associated to the
dual groups $G_c$ and $\LG_c$, with surface operators gives rise to an
equivalence of categories of $A$- and $B$-branes on ${\mc M}_H(G,x)$
and ${\mc M}_H(\LG,x)$. The mirror dual for a generic 0-brane on ${\mc
M}_H(\LG,x)$ is the $A$-brane consisting of a Hitchin fiber and a flat
unitary line bundle on it, as in the unramified case.

The analysis of \cite{GW} leads to many of the same conclusions as
those obtained in \cite{FG} by using representations of affine
Kac--Moody algebras and two-dimensional conformal field theory.

Gukov and Witten also considered \cite{GW1} more general surface
operators associated to coadjoint orbits in $\g$ and $^L\g$. The
$S$-duality between these surface operators leads to some non-trivial
and unexpected relations between these orbits. Gukov and Witten
present many interesting examples of this in \cite{GW1} drawing
connections with earlier work done by mathematicians.

In \cite{Witten:wild}, Witten has generalized the analysis of
\cite{GW,GW1} to the case of wild ramification.

\section{MORE GENERAL BRANES}    \label{BRANES}

In the previous section we discussed applications of Mirror Symmetry
of the dual Hitchin fibrations to the geometric Langlands
correspondence. We saw that the $A$-branes associated to the
$B$-branes supported at the generic flat $\LG$-bundles have very
simple description: these are the Hitchin fibers equipped with flat
unitary line bundles. But what about the $B$-branes supported at more
general flat $\LG$-bundles? Can we describe explicitly the $A$-branes
dual to them?

This question goes to the heart of the subtle interplay between
physics and mathematics of Langlands duality. Trying to answer this
question, we will see the limitations of the above analysis and a way
for its generalization incorporating more general branes. This will
lead us to surprising physical interpretation of deep mathematical
concepts such as endoscopy and Arthur's $SL_2$.

The generic flat $\LG$-bundles, for which Mirror Symmetry works so
nicely, are the ones that have no automorphisms (apart from those
coming from the center of $\LG$). They correspond to smooth points of
${\mc Y}(\LG)$ such that the corresponding Hitchin fiber is also
smooth. We should consider next the singularities of ${\mc
Y}(\LG)$. The simplest of those are the orbifold singularities. We will
discuss them, and their connection to endoscopy, in the next
subsection, following \cite{FW}. We will then talk about more general
singularities corresponding to flat $\LG$-bundles with continuous
groups of automorphisms, and what we can learn about the corresponding
categories from physics.

\subsection{Geometric Endoscopy}    \label{endoscopy}

We start with the mildest possible singularities in ${\mc Y}(\LG)$;
namely, the orbifold singular points.  The corresponding flat
$\LG$-bundles are those having finite groups of automorphisms (modulo
the center). In the classical Langlands correspondence the analogous
Galois representations are called endoscopic. They, and the
corresponding automorphic representations, play an important role in
the stabilization of the trace formula.

The simplest example, analyzed in \cite{FW} and dubbed ``geometric
endoscopy'', arises when $\LG=PGL_2$, which contains $O_2 = \Z_2
\ltimes \C^\times$ as a subgroup. Suppose that a flat $PGL_2$-bundle
$\E$ on our curve $X$ is reduced to this subgroup. Then generically it
will have the group of automorphisms $\Z_2 = \{ 1,-1 \} \subset
\C^\times$, which is the center of $O_2$ (note that the center of
$PGL_2$ itself is trivial). Therefore the corresponding points of
${\mc Y}(\LG)$ are $\Z_2$-orbifold points. This means that the
category of $B$-branes supported at such a point is equivalent to the
category $\Rep(\Z_2)$ of representations of $\Z_2$.\footnote{The
corresponding derived category has a more complicated structure, but
we will not discuss it here.} Thus, it has two irreducible objects.
Therefore we expect that the dual category of $A$-branes should also
have two irreducible objects. In fact, it was shown in \cite{FW} that
the dual Hitchin fiber has two irreducible components in this case,
and the sought-after $A$-branes are {\em fractional branes} supported
on these two components.

This was analyzed very explicitly in \cite{FW} in the case when $X$ is
an elliptic curve. Here we allow a single point of tame ramification
(along the lines of \secref{ram}) -- this turns out to be better for
our purposes than the unramified case. The corresponding Hitchin
moduli spaces are two-dimensional. They fiber over the same
one-dimensional vector space, and the fibers over all but three points
in the base are smooth elliptic curves. The three pairs of dual
singular fibers look as follows:

\vspace*{-25mm}

\begin{center}
\setlength{\unitlength}{0.6mm}
\begin{picture}(30,30)(-70,90)\label{kart}
\allinethickness{1.5pt}
\put(-24.5,35){\circle*{2}}     
\put(-5,35){\circle{40}}
\put(5,35){\circle{60}}

\thinlines
\qbezier[15](15,35)(28,40)(35,35)
\qbezier(15,35)(22,30)(35,35)
\qbezier[15](9,49)(12,58)(20,60)
\qbezier(9,49)(18,52)(20,60)
\qbezier[15](9.5,20.5)(18,18)(20,10)
\qbezier(9.5,20.5)(12,12)(20,10)

\end{picture}


\begin{picture}(90,90)(40,0)
\linethickness{1pt}
\qbezier(45,0)(-60,35)(45,70)
\qbezier(45,0)(10,35)(45,70)
\qbezier(45,0)(150,35)(45,70)
\qbezier(45,0)(80,35)(45,70)
\thinlines
\qbezier[30](-7.5,35)(15,45)(27.5,35)
\qbezier(-7.5,35)(5,25)(27.5,35)
\qbezier[20](10,15)(25,20)(33.5,15)
\qbezier(10,15)(18,10)(33.5,15)
\qbezier[20](10,55)(25,60)(33.5,55)
\qbezier(10,55)(18,50)(33.5,55)
\qbezier[30](62.5,35)(85,45)(97.5,35)
\qbezier(62.5,35)(75,25)(97.5,35)
\qbezier[20](57,15)(72,20)(80,15)
\qbezier(57,15)(65,10)(80,15)
\qbezier[20](57,55)(72,60)(80,55)
\qbezier(57,55)(65,50)(80,55)
\end{picture}
\end{center}

\noindent
\qquad \hspace*{25mm}{\small Singular Hitchin fiber in} \hspace*{24mm}
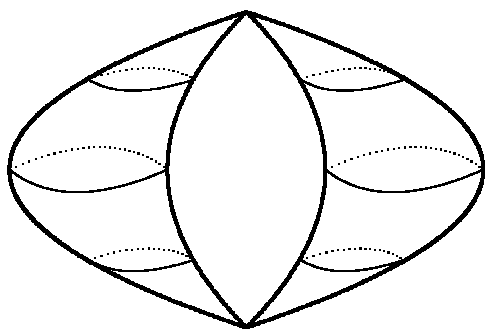
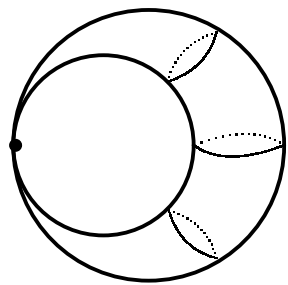
{\small Singular Hitchin fiber in}

\noindent
\qquad \hspace*{25mm}{\small the $A$-model, $G=SL_2$.} \hspace*{25mm}
  {\small the $B$-model, $^L\neg G=SO_3$.}

\vspace*{5mm}

The fiber on the $B$-model side is a projective line with a double
point, corresponding to a flat $PGL_2$-bundle that is reduced to the
subgroup $O_2$. It is a $\Z_2$-orbifold point of the moduli space. The
dual fiber on the $A$-model side is the union of two projective lines
connected at two points. These two singular points of the fiber are
actually smooth points of the ambient moduli space.

There are two irreducible $B$-branes supported at each of the
$\Z_2$-orbifold points, corresponding to two irreducible
representations of $\Z_2$. Let us denote them by $\B_+$ and
$\B_-$. The corresponding fiber of the Hitchin fibration for $SL_2$ is
the union of two components $\FF_1$ and $\FF_2$, and accordingly in
the dual $A$-model there are two irreducible $A$-branes, $\A_1$ and
$\A_2$ supported on these components (each component is a copy of
$\pone$, and therefore the only flat unitary line bundle on it is the
trivial one). These $A$-branes are dual to the $B$-branes $\B_+$ and
$\B_-$. Unlike $\B_+$ and $\B_-$, they are indistinguishable. An
apparent contradiction is explained by the fact that in the
equivalence \eqref{branes} of the categories of $A$-branes and
$B$-branes there is a twist by a $\Z_2$-gerbe which is not canonically
trivialized (see Section 9 of \cite{FW} and the footnotes on pages
\pageref{foot} and \pageref{foot1}). In order to set up an
equivalence, we need to pick a trivialization of this gerbe, and this
breaks the symmetry between $\A_1$ and $\A_2$. We also have a similar
picture when $X$ has higher genus (see \cite{FW}).

What happens when we act on $\cal B_+$ or $\cal B_-$ by the Wilson
operator $W_x, x \in X$, corresponding to the three-dimensional
adjoint representation of $\LG=PGL_2$?  Since $\cal B_+$ and $\cal
B_-$ both have skyscraper support at the same point $\E$ of ${\mc
Y}(\LG)$, $W_x$ acts on either of them by tensor product with the
three-dimensional vector space $\E_x$, the fiber of $\E$ at $x$ (in
the adjoint representation).  However, we should be more precise to
keep track of the $\Z_2$-action. Recall that the structure group of
our flat $PGL_2$-bundle $\E$ is reduced to the subgroup $O_2 = \Z_2
\ltimes \C^\times$.  Denote by $U$ the defining two-dimensional
representation of $O_2$. Then $\det U$ is the one-dimensional sign
representation induced by the homomorphism $O_2 \to \Z_2$. The adjoint
representation of $PGL_2$ decomposes into the direct sum
$$
(\det U \otimes I) \oplus (U \otimes S)
$$
as a representation of $O_2 \times \Z_2$, where $\Z_2$ is the
centralizer of $O_2$ in $PGL_2$ (the center of $O_2$), $S$ is the
sign representation of $\Z_2$, and $I$ is the trivial representation
of $\Z_2$. Therefore we have the following decomposition of the
corresponding flat vector bundle:
\begin{equation}    \label{WE}
(\det U_\E \otimes I) \oplus (U_\E \otimes S),
\end{equation}
and there is a decomposition $U_{\E}|_x\oplus \det\,U_{\E}|_x$, where
the non-trivial element of $\Z_2$ acts as $-1$ on the first summand
and as $+1$ on the second summand.  So we have
\begin{equation}    \label{elko}
W_x\cdot \cal B_\pm = \left( \cal
  B_\mp\otimes U_{\E}|_x \right) \oplus \left( \cal B_\pm\otimes \det
  U_{\E}|_x \right).
\end{equation}
Thus, individual branes $\B_+$ and $\B_-$ are not eigenbranes of the
Wilson operators; only their sum $\B_+ \oplus \B_-$ (corresponding to
the regular representation of $\Z_2$) is.

The mirror dual statement, shown in \cite{FW}, is that we have a
similar formula for the action of the corresponding 't Hooft operators
on the $\A$-branes $\A_1$ and $\A_2$. Again, only their sum (or union),
which gives the entire Hitchin fiber, is an eigenbrane of the 't Hoof
operators.

Since an eigenbrane ${\mc A}$ decomposes into two irreducible branes
${\mc A}_1$ and ${\mc A}_2$, the corresponding Hecke eigensheaf ${\mc
F}$ on $\Bun_G$ should also decompose as a direct sum of two ${\mc
D}$-modules, ${\mc F}_1$ and ${\mc F}_2$, corresponding to ${\mc A}_1$
and ${\mc A}_2$, respectively. Furthermore, these two ${\mc
D}$-modules should then separately satisfy an analogue of formula
\eqref{elko}, which is a natural modification of the standard Hecke
property. We called it in \cite{FW} the {\em fractional Hecke
property}, and the ${\mc D}$-modules ${\mc F}_1$ and ${\mc F}_2$ {\em
fractional Hecke eigensheaves}. We have also generalized this notion
to other groups in \cite{FW}.

Thus, the Mirror Symmetry picture leads us to predict the existence of
fractional Hecke eigensheaves for endoscopic flat
$\LG$-bundles.\footnote{In the case of $SL_2$ (as well as $GSp_4$) the
existence of these ${\mc D}$-modules follows from the work of Lysenko
\cite{Ly}, but for other groups this is still a conjecture.} This has
non-trivial consequences even for curves over $\Fq$, some of which
have been verified in \cite{FW}. In addition, we obtain a relation
between the group $\pi_0(\PP_b)$ of components of the generalized Prym
variety $\PP_b$, which is an open dense part of the singular Hitchin
fiber $\FF_b$ arising in the $A$-model, and the group of automorphisms
of the endoscopic flat $\LG$-bundles which are the singular points in
the dual Hitchin fiber $^L\FF_b$ of the $B$-model (this relation was
independently observed by B.C. Ng\^o). Roughly speaking, elements of
$\pi_0(\PP_b)$ label the components of $\FF_b$, and hence fractional
$A$-branes. The $B$-branes dual to them correspond to characters of
the group of automorphisms of an endoscopic flat $\LG$-bundle,
viewed as a point in $^L\FF_b$. Hence $\pi_0(\PP_b)$ should be dual
(as an abelian group) to this group of automorphisms.

The upshot of all this is that by analyzing the categories of
$A$-branes supported on the singular Hitchin fibers, we learn many
things about the geometric Langlands correspondence (and even the
Langlands correspondence for curves over finite fields) which would
have been very difficult to see directly using the conventional
formalism of ${\mc D}$-modules. This is a good illustration of the
power of this new method.

There is a link between our analysis and the classical theory of
endoscopy, due to the fact that the geometry we use is similar to that
exploited by B.C. Ng\^o in his recent proof of the fundamental lemma
\cite{Ngo}. Ng\^o has discovered a striking connection between
the orbital integrals appearing on the geometric side of the trace
formula and the ($\ell$-adic) cohomology of the Hitchin fibers in the
moduli space ${\mc M}_H(G)$, but for curves over $\Fq$; more
specifically, its decomposition under the action of the group
$\pi_0(\PP_b)$.\footnote{More precisely, Ng\^o considers a
generalization of ${\mc M}_H(G)$ parametrizing meromorphic Higgs
fields with a sufficiently large divisor of poles.}

However, there are important differences. First of all, we work over
$\C$, whereas Ng\^o works over $\Fq$. In the latter setting there is
no obvious analogue of the Homological Mirror Symmetry between ${\mc
M}_H(G)$ and ${\mc M}_H(\LG)$. Second, and more importantly, the
objects we assign to the connected components of the singular Hitchin
fiber $\FF_b$ -- the $A$-branes -- are objects of automorphic nature;
we hope to relate them to Hecke eigensheaves and ultimately to the
automorphic functions in the classical theory. Thus, these objects
should live on the {\em spectral side} of the trace formula. On the
other hand, in Ng\^o's work Hitchin fibers appear on the {\em
geometric side} of the trace formula (more precisely, its Lie algebra
version), through orbital integrals.

This raises the following question: could there be a more direct link
between individual Hitchin fibers in the moduli space $\MH(G)$ over
$\Fq$ and automorphic representations? In other words, could it be
that the passage from $A$-branes to Hecke eigensheaves discussed above
has an analogue in the classical theory as a passage from orbital
integrals to Hecke eigenfunctions? If so, then the Mirror Symmetry
picture would give us valuable insights into the Langlands
correspondence.

\subsection{$S$-duality of more general boundary conditions}
\label{more general}

More general flat $\LG$-bundles have continuous groups of
automorphisms. For instance, generic flat bundles reduced to a Cartan
subalgebra $\LH$ have the group of automorphisms $\LH$. Or consider
the trivial flat $\LG$-bundle, whose group of automorphisms is $\LG$
itself. What are the $A$-branes corresponding to these flat
$\LG$-bundles?

The picture of two dual Hitchin fibrations discussed above is too
naive to answer this question. The reason is that even if the flat
bundles with continuous groups of automorphisms are semi-stable (which
is not necessarily the case), they correspond to points of ${\mc
Y}(\LG)$ with singularity so severe that the category of $B$-branes
corresponding to it cannot be described solely in terms of the moduli
space ${\mc Y}(\LG)$. In fact, the definition of the sigma model
itself is problematic for singular target manifolds.

As an illustration, consider the quotient $\C^n/\C^\times$. The origin
has the group of automorphisms $\C^\times$. What is the category of
$B$-branes associated to this point?  Because it is a singular point,
there is no obvious answer (unlike the case of smooth points or
orbifold points, discussed above). However, we can resolve the
singularity by blowing it up. On general ground one can argue that
this resolution will not change the category of $B$-branes. The
category of $B$-branes after the resolution of singularities is the
category of coherent sheaves on ${\mathbb P}^{n-1}$. So the singular
point in $\C^n/\C^\times$ has ``swallowed'' an entire projective
space! Likewise, singular points in ${\mc Y}(\LG)$ also have
complicated ``inner structure'' which needs to be uncovered to do
justice to the corresponding categories of $B$-branes.

In order to understand better what is going on we should go back to
the four-dimensional gauge theory and look more closely at the
$S$-duality of boundary conditions there.  From the physics
perspective, this is the ``master duality'' and everything should
follow from it. The Mirror Symmetry of the Hitchin fibrations is but
the first approximation to the $S$-duality when we compactify the
theory to two dimensions.

It is instructive to recall how one obtains the Hitchin moduli spaces
in the first place: Each of the $S$-dual gauge theories has a
differential $Q$ such that $Q^2=0$, and we study the corresponding
topological field theories. In the topological theory the path
integral localizes on the moduli space of solutions to the ``BPS
equations'', which read $Q \cdot \Psi = 0$, for all fermionic fields
$\Psi$ of our theory (since $Q$ is fermionic and we want the equations
on the bosonic degrees of freedom). After that we make dimensional
reduction of these equations. This means that we assume that our
four-manifold has the form $\Sigma \times X$ and the fields on
$\Sigma$ vary ``slowly'' along $\Sigma$. The corresponding equations
have been written in \cite{KW}:
\begin{align}    \label{BPS}
F_A - \phi \wedge \phi &= 0, \\ \notag
d_A \phi = d_A \star \phi &= 0,
\end{align}
where $d_A$ is the exterior derivative corresponding to the connection
$A$, and $\star$ is the Hodge star operator. These are precisely the
Hitchin equations \cite{Hit} describing the moduli spaces $\MH(G)$ or
$\MH(\LG)$ (depending on which side of $S$-duality we are on). For
example, points of $\MH(\LG)$ in the complex structure $J$ are
semi-stable flat $\LG$-bundles on $X$. The flat connection on this
bundle is given by the formula $\nabla = A + i \phi$ (the flatness of
$\nabla$ is a corollary of \eqref{BPS}). This is how the ($B$-twisted)
sigma model on $\Sigma$ with values in $\MH(\LG)$ appears in this
story. One obtains the ($A$-twisted) sigma model with target $\MH(G)$
in the complex structure $I$ in a similar way.

However, as Kapustin and Witten explain in \cite{KW}, this derivation
breaks down when we encounter singularities of the Hitchin moduli
spaces. Thus, the sigma models with the targets $\MH(\LG)$ and
$\MH(G)$ are only approximations to the true physical theory. To
understand what happens at the singularities we have to go back to the
four-dimensional theory and analyze it more carefully (for more on
this, see \cite{Witten:rev}).

There is also another problem: in the above derivation we have not
taken into account all the fields of the super--Yang--Mills theory. In
fact, there are additional scalar fields, denoted by $\sigma$ and
$\ol\sigma$ in \cite{KW}, which we have ignored so far. The field
$\sigma$ is a section of the adjoint bundle $\g_E$ on $X$, and
$\ol\sigma$ is its complex conjugate. On the $B$-model side, which we
have so far approximated by the sigma model with the target
$\MH(\LG)$, we obtain from the BPS equations that $\sigma$ is
annihilated by the flat connection $\nabla = A + i \phi$, that is,
$\nabla \cdot \sigma = 0$. In other words, $\sigma$ belongs to the Lie
algebra of infinitesimal automorphisms of the flat bundle
$(E,\nabla)$.

Up to now we have considered generic flat $\LG$-bundles which have no
non-trivial infinitesimal automorphisms. For such flat bundles we
therefore have $\sigma \equiv 0$, and so we could safely ignore it. But
for flat bundles with continuous automorphisms this field starts
playing an important role.

The upshot of this discussion is that when we consider most general
flat bundles there are new degrees of freedom that have to be taken
into account. In order to find a physical interpretation of the
geometric Langlands correspondence for such flat bundles we need to
consider the $S$-duality of boundary conditions in the
four-dimensional gauge theory with these degrees of freedom included.

A detailed study of $S$-duality of these boundary conditions has been
undertaken by Gaiotto and Witten \cite{GaW1,GaW2}. We will only
mention two important aspects of this work.

First of all, Gaiotto and Witten show that in the non-abelian gauge
theory the $S$-dual of the Neumann boundary condition is not the usual
Dirichlet boundary condition as one might naively hope, but a more
complicated boundary condition in which the field $\sigma$ has a pole at
the boundary. This boundary condition corresponds to a solution of the
Nahm equations, which is in turn determined by an embedding of the Lie
algebra $\sw_2$ into $^L\g$ (see \cite{Witten:atiyah}).

This is parallel to the appearance of Arthur's $SL_2$ in the classical
Langlands correspondence. Arthur has conjectured that the true
parameters for ($L$-packets of) unitary automorphic representations of
$G(\AD)$ are not homomorphisms $\on{Gal}(\ol{F}/F) \to \LG$, but
rather $\on{Gal}(\ol{F}/F) \times SL_2 \to \LG$. The homomorphisms
whose restriction to the $SL_2$ factor are trivial should correspond
to the so-called tempered representations. (In the case of $GL_n$ all
cuspidal unitary representations are tempered, and that is why
Arthur's $SL_2$ does not appear in the theorem of Drinfeld and
Lafforgue quoted in \secref{LP}.) An example of non-tempered unitary
representation is the trivial representation of $G(\AD)$. According to
\cite{Arthur}, the corresponding parameter is the homomorphism
$\on{Gal}(\ol{F}/F) \times SL_2 \to \LG$, which is trivial on the
first factor and is the principal embedding of the $SL_2$ factor.

In the geometric Langlands correspondence, Arthur's $SL_2$ may be
observed in the following way. The analogue of the trivial
representation is the constant sheaf ${\mb C}$ on $\Bun_G$. It is a
Hecke eigensheaf, but the ``eigenvalues'' are complexes of vector
spaces with cohomological grading coming from the Cartan subalgebra of
the principal $SL_2$ in $\LG$. For example, consider the case of
$G=GL_n$ and let us apply the Hecke operators $H_{\check\omega_1,x}$
defined by formula \eqref{formula H1} to the constant sheaf. It
follows from the definition that
$$
H_{\check\omega_1,x}({\mb C}) \simeq H^\bullet({\mathbb P}^{n-1},\C)
\otimes {\mb C}.
$$
Thus, the eigenvalue is a graded $n$-dimensional vector space with
one-dimensional pieces in cohomological degrees
$0,2,4,\ldots,2(n-1)$. In the standard normalization of the Hecke
operator the cohomological grading is shifted to
$-(n-1),-(n-3),\ldots,(n-1)$ -- as in the grading on the
$n$-dimensional representation of $GL_n$ coming from the principal
$SL_2$.

The $A$-brane corresponding to the constant sheaf on $\Bun_G$ is the
Lagrangian submanifold of $\MH(G)$ defined by the equation $\phi=0$
(this is the zero section of the cotangent bundle to the moduli space
of semi-stable $G$-bundles inside $\MH(G)$). This $A$-brane
corresponds to a Neumann boundary condition in the 4D gauge theory
with gauge group $G_c$. According to Gaiotto and Witten, the dual
boundary condition (in the theory with gauge group $\LG_c$) is a
generalization of the Dirichlet boundary condition, in which the field
$\sigma$ has a pole at the boundary solving the Nahm equations
corresponding to the principal $SL_2$ embedding into $\LG$. Thus, we
obtain a beautiful interpretation of Arthur's $SL_2$ from the point of
view of $S$-duality of boundary conditions in gauge theory. For more
on this, see \cite{Witten:atiyah,FGukov}.

Another important feature discovered in \cite{GaW1,GaW2} is that the
$S$-duals of the general Dirichlet boundary conditions involve
coupling the 4D super--Yang-Mills to 3D superconformal QFTs at the
boundary. This means that there are some additional degrees of freedom
that we have to include to describe the geometric Langlands
correspondence.

What we learn from all this is that the true moduli spaces arising in
the $S$-duality picture are not the Hitchin moduli spaces ${\mc
M}_H(G)$ and ${\mc M}_H({}\LG)$, but some enhanced versions $\wt{\mc
M}_H(G)$ and $\wt{\mc M}_H({}\LG)$ thereof, including, in addition to
the Higgs bundle $(E,\phi)$, an element $\sigma$ in the Lie algebra of
its infinitesimal automorphisms as well as other data. (This will be
discussed in more detail in the forthcoming paper
\cite{FGukov}.) Physically, the field $\sigma$ has non-zero ``ghost
number'' $2$. Mathematically, this means that these additional degrees
of freedom have cohomological grading $2$, and so $\wt{\mc
M}_H(G)$ and $\wt{\mc M}_H({}\LG)$ are actually differential graded
(DG) stacks. Similar DG stacks have been recently studied in the
context of the categorical Langlands correspondence by V. Lafforgue
\cite{VLaff}.

\bigskip

Thus, $S$-duality of super--Yang--Mills theory offers new insights into
the Langlands correspondence and surprising new connections to
geometry. Ultimately, we have to tackle the biggest question of all:
what is the underlying reason for the Langlands duality? On the
physics side the corresponding question is: why $S$-duality? In fact,
physicists have the following elegant explanation (see
\cite{Witten:6dim}): there is a mysterious six-dimensional quantum
field theory which gives rise to the four-dimensional
super--Yang--Mills upon compactification on an elliptic curve. Roughly
speaking (the argument should be modified slightly for non-simply
laced groups), this elliptic curve is $E = \C/(\Z + \Z \tau)$, where
$\tau$ is the coupling constant of this Yang--Mills theory. Since the
action of $SL_2(\Z)$ on $\tau$ does not change the elliptic curve $E$,
we obtain that there are equivalences between the super--Yang--Mills
theories whose coupling constants are related by the action of
$SL_2(\Z)$. This should explain the $S$-duality, which corresponds to
the transformation $\tau \mapsto -1/\tau$, and hence the geometric
Langlands correspondence. But that's a topic for a future S\'eminaire
Bourbaki.

\end{document}